\documentclass[11pt,leqno]{article}
\usepackage{amsmath,amsthm}
\theoremstyle{plain}
\newtheorem{thm}{Theorem}[section]
\newtheorem{prop}[thm]{Proposition}
\newtheorem{lem}[thm]{Lemma}
\newtheorem{cor}[thm]{Corollary}
\theoremstyle{definition}
\newtheorem{defn}[thm]{Definition}
\newtheorem{exmp}[thm]{Example}
\theoremstyle{remark}
\newtheorem{rem}[thm]{Remark}
\numberwithin{equation}{section}

\renewcommand{\eqref}[1]{(\ref{#1})}
\newcommand{\xiluv}{\xi_\lambda(u)(v)}
\newcommand{\xilmpi}{\xi^{(L, M, \pi)}}
\newcommand{\xilmpidash}{\xi^{(L', M', \pi')}}
\newcommand{\etalvu}{\eta_\lambda(v)(u)}
\newcommand{\etalmpi}{\eta^{(L, M, \pi)}}
\newcommand{\etalmpidash}{\eta^{(L', M', \pi')}}
\newcommand{\sigmalmpi}{{\sigma^{(L, M, \pi)}}}
\newcommand{\Rlmpi}{{R^{(L, M, \pi)}}}
\newcommand{\Rlvmpi}{{R^{(L_V, M, \pi)}}}
\newcommand{\Rlmpidash}{{R^{(L', M', \pi')}}}
\newcommand{\RV}{{R^V}}
\newcommand{\Perm}[1]{P_{#1}}
\newcommand{\cA}{{\mathcal{A}}}
\newcommand{\cD}{{\mathcal{D}}}
\newcommand{\Hom}{{\mathrm{Hom}}}
\newcommand{\id}{{\mathrm{id}}}

\newcommand{\bsl}[1]{\backslash_{#1}\,}
\newcommand{\bsG}{{\bsl{G}}}

\newcommand{\bsL}{{\bsl{L}}}

\newcommand{\bsLV}{{\bsl{L_V}}}

\newcommand{\xiu}[2]{{\ensuremath{\xi_{#1}({#2})}}}
\newcommand{\xilu}{{\ensuremath{\xi_{\lambda}(u)}}}
\newcommand{\etav}[2]{{\ensuremath{\eta_{#1}({#2})}}}
\newcommand{\etalv}{{\ensuremath{\eta_{\lambda}(v)}}}

\begin{document}
\begin{flushleft}
\end{flushleft}
\begin{center}
{\Large
Dynamical Yang-Baxter Maps\\[.5cm]
with
an Invariance Condition}
\\[1cm]
By
\\[1cm]
Youichi Shibukawa\footnote[1]{%
Department of Mathematics, Faculty of Science,
Hokkaido University,
Sapporo 060-0810, Japan
\endgraf
e-mail: shibu@math.sci.hokudai.ac.jp}
\end{center}
\begin{abstract}
By means of left quasigroups $L=(L, \cdot)$
and ternary systems,
we construct dynamical Yang-Baxter maps
associated with $L$, $L$, and $(\cdot)$
satisfying an invariance condition that
the binary operation $(\cdot)$ of the left quasigroup $L$ defines.
Conversely, this construction characterize
such dynamical Yang-Baxter maps.
The unitary condition
of the dynamical Yang-Baxter map is discussed.
Moreover, we establish
a correspondence between two dynamical
Yang-Baxter maps constructed in this paper.
This correspondence produces a version of the
vertex-IRF correspondence.
\footnote[0]{This study was partially supported
by the Ministry
of Education, Science, Sports and Culture, Japan,
grant-in-aid for Young Scientists (B), 15740001, 2005.}
\end{abstract}
\section{Introduction}
\label{sect:introduction}
Much attention has been directed to the quantum dynamical Yang-Baxter equation (QDYBE),
a generalization of the quantum Yang-Baxter equation (QYBE)
(for example, see \cite{etingofbook}).
The dynamical Yang-Baxter map (dynamical YB map)
\cite{shibukawa05} is a set-theoretical solution to a version of the QDYBE.

Let $H$ and $X$ be nonempty sets, and
$\phi$ a map from $H\times X$ to $H$.
A map $R(\lambda):X\times X\to X\times X$
$(\lambda\in H)$
is a {\em dynamical YB map\/}
associated with $H$, $X$, and $\phi$,
iff, for every $\lambda\in H$,
$R(\lambda)$ satisfies the following equation on
$X\times X\times X$.
\begin{equation}\label{eqn:introduction:defdYB}
R_{23}(\lambda)R_{13}(\phi(\lambda, X^{(2)}))
R_{12}(\lambda)
=
R_{12}(\phi(\lambda, X^{(3)}))R_{13}(\lambda)
R_{23}(\phi(\lambda, X^{(1)})).
\end{equation}
Here
$R_{12}(\lambda)$, $R_{12}(\phi(\lambda, X^{(3)}))$,
$R_{23}(\phi(\lambda, X^{(1)}))$,
and others are the maps from $X\times X\times X$ to
itself defined
as follows:
for $u, v, w\in X$,
{\allowdisplaybreaks
\begin{eqnarray*}
&&
R_{12}(\lambda)(u, v, w)
=
(R(\lambda)(u, v), w);
\\
&&
R_{12}(\phi(\lambda, X^{(3)}))(u, v, w)
=
R_{12}(\phi(\lambda, w))(u, v, w);
\\
&&
R_{23}(\phi(\lambda, X^{(1)}))(u, v, w)
=
(u, R(\phi(\lambda, u))(v, w)).
\end{eqnarray*}}%
If a map $R(\lambda)$ is a dynamical YB map associated with
$H$, $X$, and $\phi$,
then we denote it
by $(R(\lambda); H, X, \phi)$.
Two dynamical YB maps $(R^{(1)}(\lambda_1); H_1, X_1, \phi_1)$
and
$(R^{(2)}(\lambda_2); H_2, X_2, \phi_2)$
are equal,
iff
\begin{equation}\label{eqn:intro:defequal}
H_1=H_2,\ X_1=X_2,\ \phi_1=\phi_2,\ 
\mbox{and}\ 
R^{(1)}(\lambda)=R^{(2)}(\lambda)
\end{equation}
for all $\lambda\in H_1(=H_2)$.

By the definition
of the dynamical YB map,
the Yang-Baxter map (YB map) \cite{veselov03},
a set-theoretical solution to the QYBE \cite{drinfeld92,weinstein92},
is a dynamical YB map that is independent of the
dynamical parameter $\lambda$.
Geometric crystals \cite{etingof03},
crystals \cite{hatayama02},
and semigroups of I-type
\cite{gateva04}
produce YB maps,
and so do bijective 1-cocycles
\cite{etingof99,lu00}.
Let $A$ and $G$ be groups such that $A$ acts on $G$,
and $\pi: A\to G$ a bijective 1-cocycle of the group $A$
with coefficients in the group $G$.
This triplet $(A, G, \pi)$ gives birth to a bijective YB map
\cite[Theorems $1$ and $2$]{lu00}
with the invariance condition
\eqref{eqn:characterization:invarianceybmap}
(this invariance condition is called the compatibility condition in \cite{lu00}).

By generalizing this method,
dynamical YB maps were constructed in \cite{shibukawa05}.
Let $LP=(LP, \cdot, e_{LP})$ be a loop
(see Definition \ref{defn:summary:loops}), $G=(G, *, e_G)$ a group,
and $\pi: LP\to G$ a set-theoretical bijection
satisfying $\pi(e_{LP})=e_G$.
Here $e_{LP}$ and $e_G$ are the unit elements of $LP$ and $G$,
respectively.
This triplet $(LP, G, \pi)$ produces a bijective dynamical YB map
$R^{(G)}(\lambda)$ \eqref{eqn:summary:RLP}
with the invariance condition \eqref{eqn:summary:objD3}
(see below Theorem \ref{thm:summary:3.7}).
We characterized this dynamical YB map
(Theorem \ref{thm:summary:3.7}).

However, this characterization is inadequate;
some dynamical YB maps with the
invariance condition are not constructed in this way.

This paper clarifies a characterization of dynamical YB maps with the
invariance condition.
Let $L=(L, \cdot)$ be a left quasigroup
(see Definition \ref{defn:construction:leftquasigroup}),
$M=(M, \mu)$ a ternary system
(Definition \ref{defn:construction:ternarysys})
satisfying \eqref{eqn:construction:defM1}
and
\eqref{eqn:construction:defM2},
and $\pi$ a bijection from $L$ to $M$.
The triplet $(L, M, \pi)$ produces a dynamical YB
map
$(\Rlmpi(\lambda); L, L, (\cdot))$
\eqref{eqn:construction:Rlmpi}
(Theorem \ref{thm:construction:main}).
This dynamical YB map satisfies the invariance condition
\eqref{eqn:construction:invariance}
that
the binary operation $(\cdot)$ of the left quasigroup $L$ defines.
This construction gives a characterization of
the dynamical YB maps with the invariance condition
\eqref{eqn:characterization:obD1}
(Theorem \ref{thm:characterization:main}).
If the binary operation $(\cdot)$ of $L$ is associative,
then every YB map on $L\times L$
with the invariance condition
\eqref{eqn:characterization:invarianceybmap}
is produced
by a ternary system
satisfying \eqref{eqn:construction:defM1}
and
\eqref{eqn:construction:defM2}
(see Remark \ref{rem:characterization:ybmaps}).

The organization of this paper is as follows.
After summarizing the results of the work \cite{shibukawa05}
in Section \ref{sect:summary},
we
construct the dynamical YB map
$\Rlmpi(\lambda)$
\eqref{eqn:construction:Rlmpi}
with the invariance condition 
in Section \ref{sect:construction}
(Theorem \ref{thm:construction:main}).
This dynamical YB map $\Rlmpi(\lambda)$
is a generalization of
the YB map in \cite{lu00} and the dynamical YB map in
\cite{shibukawa05}
(see Remarks \ref{rem:summary:lu} and
\ref{rem:examples:shibukawa05}).
The dynamical YB map $\Rlmpi(\lambda)$
and the corresponding dynamical braiding map
$\sigmalmpi(\lambda)$
\eqref{eqn:construction:sigmalmpi}
are expressed by means of
the 
maps $s(a)$ \eqref{eqn:construction:defsa}
and $s$ 
\eqref{eqn:construction:defs}
(see
Lemma \ref{lem:construction:sigma}
and
\eqref{eqn:unitary:Rlmpi}),
which
satisfy the braid group relation
\eqref{eqn:construction:braidrelation}.
This braid group relation
and
Lemma \ref{lem:construction:sigma}
simplify the proof
that the YB maps and the dynamical YB maps in \cite{lu00,shibukawa05}
satisfy
\eqref{eqn:introduction:defdYB}.

By means of categories $\cA$ and $\cD$
(Propositions \ref{prop:characterization:catA}
and \ref{prop:characterization:catD}),
we characterize the dynamical YB maps
associated with $L$, $L$, and $(\cdot)$
satisfying the invariance
condition \eqref{eqn:characterization:obD1}
in Sections \ref{sect:characterization} and \ref{sect:proof}
(Theorem \ref{thm:characterization:main}).

Section \ref{sect:examples}
describes several examples of the
ternary systems satisfying \eqref{eqn:construction:defM1}
and
\eqref{eqn:construction:defM2}.
For each $M=(G, \mu^G_1)$ \eqref{eqn:examples:muG1},
$(G, \mu^G_2)$ \eqref{eqn:examples:muG2},
$(G, \mu^G_3)$ \eqref{eqn:examples:muG3},
we give a characterization of the dynamical YB maps
$\Rlmpi(\lambda)$.

Sections \ref{sect:unitary}
and \ref{sect:correspondence}
deal with properties of the dynamical YB map
$\Rlmpi(\lambda)$.

Let $M$ be a ternary system constructed in Section \ref{sect:examples}.
In Section \ref{sect:unitary},
we give a necessary and sufficient condition for
the dynamical YB map $\Rlmpi(\lambda)$
to satisfy the unitary condition
\eqref{eqn:summary:unitary}
(see Propositions \ref{prop:unitary:unitary}, \ref{prop:unitary:example123},
and \ref{prop:unitary:example123-2}).
Eq.\ 
\eqref{eqn:unitary:Rlmpi}
explains the reason that
only the property
\eqref{eqn:unitary:unitary} of the
ternary system $M$ is needed in order for the dynamical YB map
$\Rlmpi(\lambda)$
to satisfy
the unitary condition
(Theorem \ref{thm:summary:unitary}).

Section \ref{sect:correspondence} gives
a correspondence between
two dynamical YB maps
called an IRF-IRF correspondence
(Proposition \ref{prop:correspondence:IRFIRF});
furthermore,
a vertex-IRF correspondence
\eqref{eqn:correspondence:vertexirf}
is discussed.
Eq.\ 
\eqref{eqn:unitary:Rlmpi}
induces this IRF-IRF correspondence.
A motivation for producing these correspondences is the exchange matrix
construction of the dynamical R-matrix
by means of the fusion matrix
(see Remark \ref{rem:correspondence:exchange}).

To end Introduction,
the author would like to thank Professor Yas-Hiro Quano
for advising him to investigate the vertex-IRF correspondence.
\section{Background information}
\label{sect:summary}
In this section, we briefly summarize the results of
the work
\cite{shibukawa05}
after introducing definitions
and notations used in this work.
\begin{defn}\label{defn:construction:leftquasigroup}
$(L, \cdot)$ is said to be a {\em left quasigroup\/}, iff $L$
is a non-empty
set, together with a binary operation $(\cdot)$
having the property that,
for all $u$, $w\in L$,
there uniquely exists
$v\in L$ such that $u\cdot v=w$
$($cf.~right quasigroups in {\rm \cite[Section I.$4$.$3$]{smith}}$)$.
\end{defn}
By this definition,
the left quasigroup $(L, \cdot)$ has another binary operation
$\bsL$
called the {\em left division\/}
\cite[Section I.$2$.$2$]{smith};
we denote by $u\bsL w$ the unique element $v\in L$
satisfying $u\cdot v=w$.
\begin{equation}\label{eqn:construction:leftdivision}
u\bsL w=v
\Leftrightarrow
u\cdot v=w.
\end{equation}
\begin{defn}\label{defn:construction:quasigroup}
A {\em quasigroup\/} $(Q, \cdot)$
is a left
quasigroup
satisfying that, for all $v$, $w\in Q$, there uniquely exists
$u\in Q$ such that $u\cdot v=w$
$($see {\rm \cite[Definition I.1.1]{pflugfelder}}
and
{\rm \cite[Section I.2]{smith}}$)$.
\end{defn}
The binary operation on a quasigroup
is not always associative.
\begin{exmp}\label{exmp:construction:123}
We define the binary operation $(*)$ on the set
$\{ 1, 2, 3\}$ of three elements by Table \ref{table:construction:three}.
Here $1*2=3$.
\begin{table}
\caption{Multiplication table of $(\{ 1, 2, 3\}, *)$
\label{table:construction:three}}
\begin{center}
\begin{tabular}{|c||c|c|c|}
\hline
$*$&$1$&$2$&$3$\\
\hline\hline
$1$&$1$&$3$&$2$\\
\hline
$2$&$2$&$1$&$3$\\
\hline
$3$&$3$&$2$&$1$\\
\hline
\end{tabular}
\end{center}
\end{table}
Then $(\{ 1, 2, 3\}, *)$ is a quasigroup, because each
element in $\{1, 2, 3\}$ appears once and only once in each row and
in each column of Table \ref{table:construction:three}
\cite[Theorem I.$1.3$]{pflugfelder}.
This binary operation $(*)$ is not associative,
since $(1*2)*3\neq 1*(2*3)$.
\end{exmp}
\begin{defn}\label{defn:summary:loops}
A {\em loop} $(LP, \cdot, e_{LP})$ is
a quasigroup $(LP, \cdot)$ satisfying that
there exists an element $e_{LP}\in LP$ such that
$u\cdot e_{LP}=e_{LP}\cdot u=u$
for all $u\in LP$
\cite[Definition I.1.10]{pflugfelder}.
\end{defn}
Because the above element $e_{LP}\in LP$
is uniquely determined,
we call $e_{LP}$ the unit element of
the loop $(LP, \cdot, e_{LP})$.

The group is a loop.
To be more precise,
the group is an associative quasigroup,
and vice versa \cite[Theorem I.$1.7$ and Definition I.$1.9$]{pflugfelder}.

We shall simply denote by $L$, $Q$, and $LP$
a left quasigroup $(L, \cdot)$,
a quasigroup $(Q, \cdot)$,
and a loop $(LP, \cdot, e_{LP})$,
respectively;
moreover,
the symbol $uv$ will be used in place of $u\cdot v$.

Next task is to demonstrate the main theorem of \cite{shibukawa05}.
Let $LP=(LP, \cdot, e_{LP})$ be a loop,
$G=(G, *, e_G)$ a group,
and $\pi:LP\to G$ a (set-theoretical) bijection
satisfying $\pi(e_{LP})=e_G$.
For $u\in LP$, we define the map $\theta(u)$
from $G$ to itself by
\begin{equation}\label{eqn:summary:theta}
\theta(u)(x)=\pi(u)^{-1}*\pi(u\pi^{-1}(x))
\quad(x\in G).
\end{equation}
Here
$\pi(u)^{-1}$ is the inverse of the element $\pi(u)$ of the group $G$.
This map $\theta(u)$ is bijective;
$\theta(u)^{-1}(x)=\pi(u\bsL\pi^{-1}(\pi(u)*x))$ $(x\in G)$.

Let $\xi_\lambda^{(G)}(u)$
and $\eta_\lambda^{(G)}(u)$
($\lambda, u\in LP$) denote
the following maps from $LP$ to itself:
for $v\in LP$,
\begin{eqnarray}
\label{eqn:summary:xishibukawa05}
&&\xi_\lambda^{(G)}(u)(v)
=\pi^{-1}\theta(\lambda)^{-1}\theta(\lambda u)\pi(v);
\\
\label{eqn:summary:etashibukawa05}
&&\eta_\lambda^{(G)}(u)(v)
=(\lambda \xi_\lambda(v)(u))\bsl{LP}((\lambda v)u).
\end{eqnarray}

Theorem $3.7$ in \cite{shibukawa05} implies
Theorem \ref{thm:summary:3.7}.
\begin{thm}
\label{thm:summary:3.7}
Let $\xilu$ and $\etalv$ $(\lambda, u, v\in LP)$
be maps from $LP$ to itself.
The following conditions are equivalent\/$:$
\begin{enumerate}
\item There exist a group $G=(G, *, e_G)$
and
a bijection $\pi: LP\to G$
satisfying
$\pi(e_{LP})=e_G$,
$\xilu=\xi_\lambda^{(G)}(u)$,
and
$\etalv=\eta_\lambda^{(G)}(v)$
for all $\lambda, u, v\in LP$\/$;$
\item
The maps $\xilu$ and $\etalv$
satisfy the properties below.
\begin{eqnarray}
&&
\xilu\xiu{\lambda u}{v}
=
\xiu{\lambda}{\lambda\bsl{LP}((\lambda u)v)}
\quad(\forall\lambda, u, v\in LP),
\notag
\\
&&
\etav{\lambda\xilu(v)}{w}(\etalv(u))
=
\etav{\lambda}{(\lambda u)\bsl{LP}(((\lambda u)v)w)}(u)
\notag
\\\notag
&&\quad\quad\quad\quad
\quad\quad\quad\quad
\quad\quad\quad\quad
\quad\quad\quad\quad
(\forall\lambda, u, v, w\in LP),
\\
&&
(\lambda \xilu(v))\etalv(u)=(\lambda u)v
\quad(\forall\lambda, u, v\in LP),
\label{eqn:summary:objD3}
\\
&&
\xiu{\lambda}{e_{LP}}=\etav{\lambda}{e_{LP}}=\mathrm{id}_{LP}
\quad(\forall\lambda\in LP).
\notag
\end{eqnarray}
\end{enumerate}
\end{thm}
We define the map $R^{(G)}(\lambda)$ $(\lambda \in LP)$
from $LP\times LP$ to itself by
\begin{equation}\label{eqn:summary:RLP}
R^{(G)}(\lambda)(u, v)=(\eta_\lambda^{(G)}(v)(u),
\xi_\lambda^{(G)}(u)(v))
\quad(u, v\in LP).
\end{equation}
From Propositions $3.3$ and $5.1$ in \cite{shibukawa05}
and the above theorem,
this map $R^{(G)}(\lambda)$
is a bijective dynamical YB map associated
with $LP$, $LP$, and $(\cdot)$.
\begin{rem}\label{rem:summary:lu}
This method produces all YB maps constructed in the work
\cite{lu00}.
We suppose that $LP$ is a group and that
the map $\theta(u)$ satisfies
\begin{equation}\label{eqn:summary:thetauv}
\theta(uv)=\theta(u)\theta(v)
\quad(\forall u, v\in LP).
\end{equation}
The definition \eqref{eqn:summary:theta}
of the map $\theta(u)$
and \eqref{eqn:summary:thetauv}
immediately induce that $\pi$ is a bijective $1$-cocycle
of $LP$ with coefficients in $G$
\cite[(8)]{lu00}.
Because of \eqref{eqn:summary:xishibukawa05}
and
\eqref{eqn:summary:thetauv},
the map $\xi_\lambda^{(G)}(u)$
\eqref{eqn:summary:xishibukawa05}
is independent of the dynamical parameter $\lambda$.
Since $LP$ is a group,
the map
$\eta_\lambda^{(G)}(v)$
\eqref{eqn:summary:etashibukawa05}
is also independent of $\lambda$.
By the definition of the YB map in
\cite[Case 2 in the proof of Theorem 2]{lu00},
the map $R^{(G)}(\lambda)$ \eqref{eqn:summary:RLP}
is the YB map
in \cite{lu00}.
\end{rem}

Next we shall show a necessary and sufficient
condition for the dynamical YB map $R^{(G)}(\lambda)$
\eqref{eqn:summary:RLP}
to satisfy the unitary condition.

Let $R(\lambda)$ be a dynamical YB map associated with $H$, $X$, and $\phi$.
This dynamical YB map $R(\lambda)$ 
is said to satisfy the {\em unitary condition\/}
\cite[Section 5]{shibukawa05},
iff
\begin{equation}\label{eqn:summary:unitary}
R(\lambda)\Perm{X}R(\lambda)=\Perm{X}
\quad(\forall \lambda\in H).
\end{equation}
Here we denote by $\Perm{X}$ the map from $X\times X$ to itself
defined by
\begin{equation}\label{eqn:construction:flip}
\Perm{X}(u, v)=(v, u)\quad(u, v\in X).
\end{equation}
\begin{thm}[Corollary $5.6$ in \cite{shibukawa05}]
\label{thm:summary:unitary}
The dynamical YB map $R^{(G)}(\lambda)$
$\eqref{eqn:summary:RLP}$
satisfies the unitary condition,
if and only if
the group $G$ is abelian.
\end{thm} 

Before ending this section,
let us introduce dynamical braiding maps
(see \cite[Section $2$]{shibukawa05}).
\begin{defn}\label{defn:summary:dbraiding}
Let $H$ and $X$ be nonempty sets, and
$\phi$ a map from $H\times X$ to $H$.
A map $\sigma(\lambda):X\times X\to X\times X$
$(\lambda\in H)$
is a
{\em dynamical braiding map\/}
associated with $H$, $X$, and $\phi$,
iff, for every $\lambda\in H$,
$\sigma(\lambda)$ satisfies the following equation
on $X\times X\times X$.
\[
\sigma(\lambda)_{12}
\sigma(\phi(\lambda, X^{(1)}))_{23}
\sigma(\lambda)_{12}
=
\sigma(\phi(\lambda,X^{(1)}))_{23}
\sigma(\lambda)_{12}
\sigma(\phi(\lambda, X^{(1)}))_{23}.
\]
\end{defn}

The concepts of the dynamical braiding map and
the dynamical YB map are exactly the same.
\begin{prop}[Proposition $2.1$ in \cite{shibukawa05}]
\label{prop:construction:2.1}
Let 
$R(\lambda)$ and $\sigma(\lambda)$
$(\lambda\in H)$ be maps from $X\times X$ to itself
satisfying
$\sigma(\lambda)=\Perm{X}R(\lambda)$
for all $\lambda\in H$.
Here $\Perm{X}$ is the map $\eqref{eqn:construction:flip}$.
The map $R(\lambda)$ is a dynamical YB
map associated with $H$, $X$, and $\phi$,
if and only if the map
$\sigma(\lambda)$ is a dynamical braiding
map associated with $H$, $X$, and $\phi$.
\end{prop}
\section{Construction}
\label{sect:construction}
Our main aim in the present section is to show how to
construct
dynamical YB maps.
This is a generalization of the works
\cite{lu00,shibukawa05}
(see Remarks \ref{rem:summary:lu}
and \ref{rem:examples:shibukawa05}). 
\begin{defn}\label{defn:construction:ternarysys}
A {\em ternary system\/} $(M, \mu)$ is a pair of a nonempty set $M$ and
a ternary operation $\mu: M\times M\times M\to M$.
\end{defn}
We shall simply denote by $M$
a ternary system $(M, \mu)$.

Let $L=(L, \cdot)$ be a left quasigroup,
$M=(M, \mu)$ a ternary system,
and $\pi:L\to M$ a (set-theoretical) bijection.
For $\lambda, u\in L$,
we define the maps $\xilmpi_\lambda(u):L\to L$
and $\etalmpi_\lambda(u):L\to L$
as follows: for $v\in L$,
\begin{eqnarray}
&&\xilmpi_\lambda(u)(v)=
\lambda\bsL\pi^{-1}(\mu(\pi(\lambda), \pi(\lambda u), \pi((\lambda u)v)));
\label{eqn:construction:xilmpi}
\\
&&
\etalmpi_\lambda(u)(v)=(\lambda\xilmpi_\lambda(v)(u))\bsL
((\lambda v)u).
\label{eqn:construction:etalmpi}
\end{eqnarray}

Let $\Rlmpi(\lambda)$ $(\lambda\in L)$ denote
the map
from $L\times L$ to itself
defined by
\begin{equation}
\Rlmpi(\lambda)(u, v)=(\etalmpi_\lambda(v)(u), \xilmpi_\lambda(u)(v))
\quad(u, v\in L).
\label{eqn:construction:Rlmpi}
\end{equation}

Since $L$ is a left quasigroup,
\eqref{eqn:construction:etalmpi}
is equivalent to
the following invariance condition
of the map $\Rlmpi(\lambda)$ (see Remark
\ref{rem:characterization:invariance}).
\begin{equation}\label{eqn:construction:invariance}
(\lambda\xilmpi_\lambda(u)(v))\etalmpi_\lambda(v)(u)
=
(\lambda u)v
\quad
(\forall \lambda, u, v\in L).
\end{equation}
\begin{thm}\label{thm:construction:main}
The map $\Rlmpi(\lambda)$
$\eqref{eqn:construction:Rlmpi}$
is a dynamical YB map associated with
$L$, $L$, and $(\cdot)$,
if and only if
the ternary system $M$ satisfies the following equations
for all $a, b, c, d\in M$\/$:$
\begin{eqnarray}
\label{eqn:construction:defM1}
&&\mu(a, \mu(a, b, c), \mu(\mu(a, b, c), c, d))
=
\mu(a, b, \mu(b, c, d));
\\
&&
\mu(\mu(a, b, c), c, d)
=
\mu(\mu(a, b, \mu(b, c, d)), \mu(b, c, d), d).
\label{eqn:construction:defM2}
\end{eqnarray}
\end{thm}
This theorem induces that the triplet $(L, M, \pi)$ with
\eqref{eqn:construction:defM1}
and
\eqref{eqn:construction:defM2}
gives birth to a dynamical YB map $\Rlmpi(\lambda)$
\eqref{eqn:construction:Rlmpi}
associated with $L$, $L$, and $(\cdot)$
satisfying the invariance condition \eqref{eqn:construction:invariance}.

Section \ref{sect:examples} describes several ternary systems with
\eqref{eqn:construction:defM1}
and
\eqref{eqn:construction:defM2}.

Let $a$ be an element of the ternary system $M$.
For the proof of Theorem \ref{thm:construction:main},
we need the maps $s(a): M\times M\to M\times M$ and
$s: M\times M\times M\to M\times M\times M$:
for $x, y, z\in M$,
\begin{eqnarray}
\label{eqn:construction:defsa}
&&s(a)(x, y)=(\mu(a, x, y), y);
\\\label{eqn:construction:defs}
&&s(x, y, z)=(x, \mu(x, y, z), z).
\end{eqnarray}
\begin{lem}\label{lem:construction:braidrelation}
The maps $s(a)_{12}$ and $s$ satisfy the braid group relation
\begin{equation}\label{eqn:construction:braidrelation}
s(a)_{12}ss(a)_{12}=ss(a)_{12}s
\quad(\forall a\in M),
\end{equation}
if and only if
the ternary system $M$ satisfies $\eqref{eqn:construction:defM1}$
and $\eqref{eqn:construction:defM2}$.
\end{lem}
\begin{proof}
The proof is straightforward.
\end{proof}
Let $\lambda$ be an element of the left quasigroup $L$,
and let $f_\lambda$ denote the following map from
$L\times L\times L$ to itself.
\begin{equation}\label{eqn:construction:flambda}
f_\lambda(u, v, w)=(\lambda u, (\lambda u)v, ((\lambda u)v)w)
\quad
(u, v, w\in L).
\end{equation}
\begin{lem}\label{lem:construction:flambda}
The map $f_\lambda$ is bijective\/$;$
$f_\lambda^{-1}(u, v, w)=(\lambda \bsL u, u\bsL v, v\bsL w)$
$(u, v, w\in L)$.
\end{lem}
We define the maps
$\sigmalmpi(\lambda): L\times L\to L\times L$
and
$\sigmalmpi(\lambda L^{(1)})_{23}: L\times L\times L\to L\times L\times L$
as follows:
for $u, v, w\in L$,
{\allowdisplaybreaks
\begin{eqnarray}
\label{eqn:construction:sigmalmpi}
&&\sigmalmpi(\lambda)(u, v)=
(\xilmpi_\lambda(u)(v),
\etalmpi_\lambda(v)(u));
\\\nonumber
&&\sigmalmpi(\lambda L^{(1)})_{23}
(u, v, w)
=
(u, \sigmalmpi(\lambda u)(v, w)).
\end{eqnarray}
}
\begin{lem}\label{lem:construction:sigma}
The maps $\sigmalmpi(\lambda)_{12}$
and $\sigmalmpi(\lambda L^{(1)})_{23}$
are expressed by means of the maps $s(\pi(\lambda))$
and $s$, respectively\/$:$
\begin{eqnarray*}
&&\sigmalmpi(\lambda)_{12}=f_\lambda^{-1}(\pi^{-1}\times\pi^{-1}\times\pi^{-1})
s(\pi(\lambda))_{12}(\pi\times\pi\times\pi)f_\lambda;
\\
&&
\sigmalmpi(\lambda L^{(1)})_{23}=f_\lambda^{-1}
(\pi^{-1}\times\pi^{-1}\times\pi^{-1})
s(\pi\times\pi\times\pi)f_\lambda.
\end{eqnarray*}
\end{lem}
\begin{proof}[Proof of Theorem $\ref{thm:construction:main}$.]
From Lemmas \ref{lem:construction:braidrelation},
\ref{lem:construction:flambda},
and \ref{lem:construction:sigma},
$\sigmalmpi(\lambda)$ is a dynamical braiding map
associated with $L$, $L$, and $(\cdot)$
(see Definition \ref{defn:summary:dbraiding}),
if and only if the ternary system $M$ satisfies
\eqref{eqn:construction:defM1}
and
\eqref{eqn:construction:defM2}.

Proposition
\ref{prop:construction:2.1}, \eqref{eqn:construction:Rlmpi},
and \eqref{eqn:construction:sigmalmpi}
complete the proof.
\end{proof}
\section{Characterization}
\label{sect:characterization}
This section clarifies a characterization of the dynamical
YB map with the invariance condition \eqref{eqn:characterization:obD1};
this map is exactly the dynamical YB map
$\Rlmpi(\lambda)$ 
\eqref{eqn:construction:Rlmpi}
constructed in the previous section.

To give a characterization,
we need categories $\cA$ and $\cD$
(cf.\ \cite[Section 3]{shibukawa05}).
For category theory, see \cite{kassel,maclane98}.
Let $L=(L, \cdot)$ be a
left quasigroup
(see Definition \ref{defn:construction:leftquasigroup}),
$M=(M, \mu)$ a ternary system
(Definition \ref{defn:construction:ternarysys})
satisfying
\eqref{eqn:construction:defM1}
and
\eqref{eqn:construction:defM2},
and $\pi:L\to M$ a bijection.
We denote by $LMB$ the set of all such triplets $(L, M, \pi)$.

Triplets $(L, (M, \mu), \pi)$ and $(L', (M', \mu'), \pi')\in LMB$
are equivalent,
iff
$L=L'$ as left quasigroups
and
the map $h:=\pi'\pi^{-1}: M\to M'$
is a homomorphism of ternary systems;
that is,
the map $h: M\to M'$ satisfies
\begin{equation}
h(\mu(a, b, c))
=
\mu'(h(a), h(b), h(c))
\quad(\forall a, b, c\in M).
\label{eqn:characterization:systemhom}
\end{equation}
This is an equivalence relation,
and we write it in the form $(L, M, \pi)\sim(L', M', \pi')$.

Let $[(L, M, \pi)]$ denote the equivalence class
to which $(L, M, \pi)\in LMB$ belongs,
$Ob(\cA)$ the class of all equivalence classes with
respect to this relation.

By the definition of the relation $\sim$, all the
left quasigroups $L$ in representatives $(L, M, \pi)$ of $V\in Ob(\cA)$
are the same.
We denote by $L_V$ the left quasigroup $L$.
\begin{defn}\label{defn:characterization:homA}
Let $V$ and $V'$ be elements of $Ob(\cA)$.
We say that $f: V\to V'$ is an element of $\Hom(\cA)$,
iff
$f: L_V\to L_{V'}$ is a homomorphism
of left quasigroups such that
$\pi' f\pi^{-1}: M\to M'$ is a homomorphism
\eqref{eqn:characterization:systemhom}
of ternary systems
for all representatives
$(L_V, M, \pi)\in V$ and
$(L_{V'}, M', \pi')\in V'$.
\end{defn}
\begin{rem}\label{rem:characterization:homcA}
On account of the definition of the equivalence relation $\sim$,
$f: V\to V'\in\Hom(\cA)$,
iff
$f: L_V\to L_{V'}$ is a homomorphism
of left quasigroups
and
there exist representatives
$(L_V, M, \pi)\in V$ and
$(L_{V'}, M', \pi')\in V'$
such that 
$\pi' f\pi^{-1}: M\to M'$ is a homomorphism
of ternary systems.
\end{rem}
\begin{prop}\label{prop:characterization:catA}
$\cA$ is a category\/$:$
its objects are the elements of $Ob(\cA)$\/$;$
its morphisms are the elements of
$\Hom(\cA)$\/$;$
the identity $\id$ and the composition $\circ$
of the category $\cA$ are defined as follows\/$:$
\begin{eqnarray}
\label{eqn:characterization:idcomp1}
&&
\mbox{for}\ V\in Ob(\cA),
\id_V(u)=u\ (u\in L_V);
\\\nonumber
&&
\mbox{for}\ 
f: V\to V', g: V'\to V''\in\Hom(\cA)\ 
(V, V', V''\in Ob(\cA)),
\\\label{eqn:characterization:idcomp2}
&&
(g\circ f)(u)=g(f(u))\quad
(u\in L_V).
\end{eqnarray}
\end{prop}

The next task is to introduce a category $\cD$.
Let $L=(L, \cdot)$ be a left quasigroup,
and $R(\lambda)$ $(\lambda\in L)$ a map from $L\times L$ to itself.
We denote by
$\xi_\lambda(u)$
and
$\eta_\lambda(v)$
$(\lambda, u, v\in L)$
the following maps from $L$ to $L$.
\begin{equation}\label{eqn:characterization:etaxi}
(\eta_\lambda(v)(u), \xi_\lambda(u)(v))=R(\lambda)(u, v).
\end{equation}

Let us suppose that this map $R(\lambda)$ is
a dynamical YB map associated with
$L$,$L$, and $(\cdot)$
satisfying the invariance condition below:
\begin{equation}\label{eqn:characterization:obD1}
(\lambda\xiluv)\etalvu=(\lambda u)v
\quad
(\forall \lambda, u, v\in L).
\end{equation}
\begin{rem}\label{rem:characterization:invariance}
To be more precise,
Eq.\ \eqref{eqn:characterization:obD1}
is the invariance condition for the corresponding
dynamical braiding
map $\sigma(\lambda)=\Perm{L}R(\lambda)$
(for the map $\Perm{L}$,
see \eqref{eqn:construction:flip}).
\end{rem}
We denote by $Ob(\cD)$ the class of all such pairs
$(L, R(\lambda))$.
\begin{defn}\label{defn:characterization:homD}
Let $V=(L, R(\lambda))$ and $V'=(L', R'(\lambda'))$
be elements of $Ob(\cD)$.
We say that $f: V\to V'$ is an element of $\Hom(\cD)$, iff
$f: L\to L'$ is a homomorphism of left quasigroups satisfying
$R'(f(\lambda))(f\times f)=(f\times f)R(\lambda)$
for all $\lambda\in L$.
\end{defn}
\begin{prop}\label{prop:characterization:catD}
$\cD$ is a category\/$:$
its objects are the elements of $Ob(\cD)$\/$;$
its morphisms are the elements of
$\Hom(\cD)$\/$;$
the definitions of the identity $\id$ and the composition $\circ$ are
similar to $\eqref{eqn:characterization:idcomp1}$
and
$\eqref{eqn:characterization:idcomp2}$.
\end{prop}

Theorem \ref{thm:characterization:main} gives
a characterization of the dynamical YB maps
with the invariance condition
\eqref{eqn:characterization:obD1}.
\begin{thm}\label{thm:characterization:main}
The category $\cA$ is isomorphic to the category $\cD$.
\end{thm}

The next section will be devoted to the proof of this theorem;
we shall explicitly construct functors
$S: \cA\to\cD$
and
$T: \cD\to\cA$
satisfying $TS=\id_\cA$
and $ST=\id_\cD$.
\begin{rem}\label{rem:characterization:ybmaps}
Theorem \ref{thm:characterization:main}
produces an application of YB maps.
Let $L$ be an associative left quasigroup
and
$R$ a YB map defined on the set $L\times L$.
We denote by $\xi(u)$ and $\eta(v)$ $(u, v\in L)$
the maps from $L$ to itself defined by
$(\eta(v)(u), \xi(u)(v))=R(u, v)$.
We suppose that these maps satisfy
the invariance condition
\begin{equation}\label{eqn:characterization:invarianceybmap}
\xi(u)(v)\eta(v)(u)=uv\quad
(\forall u, v\in L).
\end{equation}
Because the binary operation of $L$ is associative,
\eqref{eqn:characterization:invarianceybmap}
is equivalent to the invariance condition
\eqref{eqn:characterization:obD1},
and
$(L, R)$ is an object of the category $\cD$
as a result.
From Theorem \ref{thm:characterization:main}
(and its proof),
this YB map $R$ is constructed by a ternary system
(Definition \ref{defn:construction:ternarysys})
satisfying
\eqref{eqn:construction:defM1}
and
\eqref{eqn:construction:defM2}.
\end{rem}
\section{Proof of Theorem \ref{thm:characterization:main}}
\label{sect:proof}
This section presents the proof of Theorem
\ref{thm:characterization:main}.
We shall first define a functor $S:\cA\to\cD$.
\begin{lem}\label{lem:proof:equiv}
Let $((L, \cdot), (M, \mu), \pi)$ and $((L', \cdot'), (M', \mu'), \pi')$
be elements of $LMB$.
The following conditions are equivalent\/$:$
\begin{enumerate}
\item
$(L, M, \pi)\sim(L', M', \pi')$\/$;$
\item
$(\Rlmpi(\lambda); L, L, (\cdot))=
(R^{(L', M', \pi')}(\lambda); L', L', (\cdot'))$
$\eqref{eqn:intro:defequal}$$;$
that is,
$L=L'$ as left quasigroups,
and
$\Rlmpi(\lambda)=R^{(L, M', \pi')}(\lambda)$
for all $\lambda\in L(=L')$.
\end{enumerate}
\end{lem}
\begin{proof}
On account of \eqref{eqn:construction:etalmpi},
the condition $(2)$ is equivalent to the condition $(3)$ below.
\begin{enumerate}
\setcounter{enumi}{2}
\item
$L=L'$ as left quasigroups,
and
$\xilmpi_\lambda(u)=\xi_\lambda^{(L, M', \pi')}(u)$
for all $\lambda, u\in L$.
\end{enumerate}

We shall only show $(1)$ from $(3)$.
It suffices to prove that the map $h=\pi'\pi^{-1}: M\to M'$ is
a homomorphism
\eqref{eqn:characterization:systemhom}
of ternary systems.

Let $a, b$, and $c$ be elements of $M$.
We define the elements $\lambda$, $u$, and $v$
of the left quasigroup $L$ by
$\lambda=\pi^{-1}(a)$,
$u=\pi^{-1}(a)\bsL\pi^{-1}(b)$,
and
$v=\pi^{-1}(b)\bsL\pi^{-1}(c)$.
Because of \eqref{eqn:construction:xilmpi},
\begin{eqnarray*}
\pi^{-1}(\mu(a, b, c))
&=&\lambda
(\lambda\bsL
\pi^{-1}(\mu(\pi(\lambda), \pi(\lambda u),
\pi((\lambda u)v))))
\\
&=&
\lambda\xilmpi_\lambda(u)(v).
\end{eqnarray*}
It follows from the condition $(3)$
that
\[
\pi^{-1}(\mu(a, b, c))
=
\pi'^{-1}(\mu'(h(a), h(b), h(c))).
\]
Hence, the map $h$ is a homomorphism of ternary systems.
\end{proof}

Let $V=[(L_V, M, \pi)]$ be an object of the category $\cA$.
From Lemma $\ref{lem:proof:equiv}$,
we can
define the dynamical YB map $\RV(\lambda)$
associated with $L_V$, $L_V$, and $(\cdot)$,
by using the dynamical YB map $\Rlvmpi(\lambda)$
\eqref{eqn:construction:Rlmpi};
\begin{equation}
\RV(\lambda)=\Rlvmpi(\lambda).
\label{eqn:proof:RV}
\end{equation}

Let $V$ be an object of the category $\cA$.
We define $S(V)$ by $S(V)=(L_V, \RV(\lambda))$.
\begin{lem}
For $V\in Ob(\cA)$,
$S(V)$ is an object of the category $\cD$.
\end{lem}
\begin{proof}
The proof is immediate from \eqref{eqn:construction:invariance}
and Theorem \ref{thm:construction:main}.
\end{proof}
\begin{lem}
Let $V$ and $V'$ be objects of the category $\cA$.
If $f: V\to V'\in\Hom(\cA)$,
then $f$ is a morphism of the category $\cD$
whose source and target are $S(V)$ and $S(V')$, respectively.
\end{lem}
\begin{proof}
Let $(L, (M, \mu), \pi)$
and $(L', (M', \mu'), \pi')$
be representatives of $V$ and $V'$, respectively.

We shall demonstrate that
$\Rlmpidash(f(\lambda))(f\times f)=(f\times f)\Rlmpi(\lambda)$
for all $\lambda\in L$.
Let $u$ and $v$ be elements of the left quasigroup $L$.
Because the map $\pi' f\pi^{-1}: M\to M'$ is a homomorphism
\eqref{eqn:characterization:systemhom}
of ternary systems
(see Definition \ref{defn:characterization:homA}),
\begin{eqnarray}\label{eqn:proof:cAcD1}
&&
f(\pi^{-1}
(\mu(\pi(\lambda), \pi(\lambda u), \pi((\lambda u)v))))
\\\nonumber
&=&
\pi'^{-1}(\mu'(\pi'(f(\lambda)), \pi'(f(\lambda u)), \pi'(f((\lambda u)v)))).
\end{eqnarray}

Since
the map $f: L\to L'$ is a homomorphism of left quasigroups,
\eqref{eqn:construction:xilmpi} and
\eqref{eqn:proof:cAcD1}
induce that
\begin{eqnarray}\label{eqn:proof:fxilmpi}
f(\xilmpi_\lambda(u)(v))
&=&
f(\lambda)\bsL f(\pi^{-1}
(\mu(\pi(\lambda), \pi(\lambda u), \pi((\lambda u)v))))
\\\nonumber
&=&
\xilmpidash_{f(\lambda)}(f(u))(f(v)).
\end{eqnarray}

The above equation and \eqref{eqn:construction:etalmpi} lead to that
\begin{equation}\label{eqn:proof:fetalmpi}
f(\etalmpi_\lambda(v)(u))=\etalmpidash_{f(\lambda)}(f(u))(f(v)),
\end{equation}
because the map $f: L\to L'$ is a homomorphism of left quasigroups.

From \eqref{eqn:proof:fxilmpi}
and
\eqref{eqn:proof:fetalmpi},
$\Rlmpidash(f(\lambda))(f\times f)=(f\times f)\Rlmpi(\lambda)$
for all $\lambda\in L$.
Thus $f: S(V)\to S(V')$ is a morphism of the category $\cD$.
\end{proof}
For $f: V\to V'\in \Hom(\cA)$,
we define $S(f): S(V)\to S(V')\in \Hom(\cD)$
by $S(f)=f$.
\begin{prop}
$S$ is a functor from the category $\cA$ to the category $\cD$.
\end{prop}

The next task is to introduce a functor $T: \cD\to\cA$.
Let $V=(L, R(\lambda))$ be an object of the category $\cD$.
We define the maps
$\xi_\lambda(u)$
and
$\eta_\lambda(v)$
$(\lambda, u, v\in L)$
from $L$ to $L$ by 
\eqref{eqn:characterization:etaxi}.
Let $\mu_L$ denote the ternary operation on $L$ defined by
\begin{equation}\label{eqn:proof:defmuL}
\mu_L(a, b, c)=
a\xi_a(a\bsL b)(b\bsL c)
\quad(a, b, c\in L).
\end{equation}
\begin{lem}
The ternary operation $\mu_L$
$\eqref{eqn:proof:defmuL}$
satisfies
$\eqref{eqn:construction:defM1}$
and
$\eqref{eqn:construction:defM2}$.
\end{lem}
\begin{proof}
Let $\lambda$ be an element of the left quasigroup $L$.
We define the maps
$s(\lambda):L\times L\to L\times L$ and
$s: L\times L\times L\to L\times L\times L$
as follows:
\begin{eqnarray*}
&&
s(\lambda)(a, b)=(\mu_L(\lambda, a, b), b)\quad (a, b\in L);
\\
&&
s(a, b, c)=(a, \mu_L(a, b, c), c)
\quad (a, b, c\in L).
\end{eqnarray*}
On account of Lemma \ref{lem:construction:braidrelation},
it suffices to prove that the maps
$s(\lambda)_{12}$
and $s$
satisfy the braid group relation
\eqref{eqn:construction:braidrelation}.

Let $\sigma(\lambda)$ denote the map from
$L\times L$ to itself defined by
$\sigma(\lambda)=\Perm{L}R(\lambda)$.
Here $\Perm{L}$ is the map \eqref{eqn:construction:flip}.
The maps $\sigma(\lambda)_{12}$
and $\sigma(\lambda L^{(1)})_{23}$
are expressed
as follows
(cf. Lemma \ref{lem:construction:sigma}):
\[
\sigma(\lambda)_{12}=f_\lambda^{-1}
s(\lambda)_{12}f_\lambda;
\quad
\sigma(\lambda L^{(1)})_{23}=f_\lambda^{-1}
sf_\lambda.
\]
Here $f_\lambda$ is the bijection \eqref{eqn:construction:flambda}
(see Lemma \ref{lem:construction:flambda}).

Because 
$\sigma(\lambda)$ is a dynamical braiding map associated with
$L$, $L$, and $(\cdot)$
(see Definition \ref{defn:summary:dbraiding}
and Proposition \ref{prop:construction:2.1}),
the maps
$s(\lambda)_{12}$
and $s$ satisfy 
\eqref{eqn:construction:braidrelation}.
This completes the proof.
\end{proof}
\begin{cor}
The triplet $(L, (L, \mu_L), \id_L)$ is an element of the set $LMB$.
\end{cor}
Let $V=(L, R(\lambda))$ be an object of the category $\cD$.
We define $T(V)\in Ob(\cA)$ by
$
T(V)=[(L, (L, \mu_L), \id_L)]
$.
\begin{lem}
Let $V$ and $V'$ be objects of the category $\cD$.
If $f: V\to V'\in\Hom(\cD)$,
then $f$ is a morphism of the category $\cA$
whose source and target are $T(V)$ and $T(V')$, respectively.
\end{lem}
\begin{proof}
Let $(L, R(\lambda))$ and $(L', R'(\lambda'))$ denote
the objects $V$ and $V'$, respectively.
We define the maps
$\xi_\lambda(u): L\to L$,
$\eta_\lambda(v): L\to L$,
$\xi'_{\lambda'}(u'): L'\to L'$,
and
$\eta'_{\lambda'}(v'): L'\to L'$
$(\lambda, u, v\in L, \lambda', u', v'\in L')$
by
\eqref{eqn:characterization:etaxi}:
$(\eta_\lambda(v)(u), \xi_\lambda(u)(v))=R(\lambda)(u, v)$;
$(\eta'_{\lambda'}(v')(u'), \xi'_{\lambda'}(u')(v'))=R'(\lambda')(u', v')$.
Let $\mu_L$ and $\mu_{L'}$ denote the ternary operations
on $L$ and $L'$ defined by \eqref{eqn:proof:defmuL}, respectively.

We shall show that
$f: (L, \mu_L)\to (L', \mu_{L'})$ is a homomorphism
\eqref{eqn:characterization:systemhom} of ternary systems.
By Definition \ref{defn:characterization:homD},
$R'(f(\lambda))(f\times f)=(f\times f)R(\lambda)$
for all $\lambda\in L$.
As a result,
$f(\xi_\lambda(u)(v))=\xi'_{f(\lambda)}(f(u))(f(v))$
for all $\lambda, u, v\in L$.
Because $f: L\to L'$ is a homomorphism of left quasigroups,
the above equation and the definition \eqref{eqn:proof:defmuL}
of $\mu_L$ and $\mu_{L'}$ induce that
$f(\mu_L(a, b, c))=\mu_{L'}(f(a), f(b), f(c))$
for all $a, b, c\in L$.
That is,
$f: (L, \mu_L)\to (L', \mu_{L'})$ is a homomorphism of ternary systems.

Since $(L, (L, \mu_L), \id_L)\in T(V)$,
$(L', (L', \mu_{L'}), \id_{L'})\in T(V')$,
and $\id_{L'}^{-1}f\id_L: (L, \mu_L)\to (L', \mu_{L'})$
is a homomorphism of ternary systems,
Remark
\ref{rem:characterization:homcA}
gives rise to that $f: T(V)\to T(V')$ is a morphism of the category $\cA$.
\end{proof}
For $f: V\to V'\in \Hom(\cD)$,
we define $T(f): T(V)\to T(V')\in \Hom(\cA)$
by $T(f)=f$.
\begin{prop}
$T$ is a functor from the category $\cD$ to the category $\cA$.
\end{prop}
\begin{proof}[Proof of Theorem $\ref{thm:characterization:main}$.]
We shall only demonstrate that
$TS(V)=V$
for $V\in Ob(\cA)$.

Let $(L_V, (M, \mu), \pi)$ be a representative of $V$.
By the definitions,
$S(V)=(L_V, R^{(L_V, M, \pi)}(\lambda))$
and $TS(V)=[(L_V, (L_V, \mu_{L_V}), \id_{L_V})]$
(see \eqref{eqn:proof:RV} and \eqref{eqn:proof:defmuL}).

For the proof, it suffices to show that
$(L_V, (L_V, \mu_{L_V}), \id_{L_V})\sim(L_V, M, \pi)$.
Let $a$, $b$, and $c$ be elements of $L_V$.
From \eqref{eqn:construction:xilmpi}
and
\eqref{eqn:proof:defmuL},
\[
\pi(\mu_{L_V}(a, b, c))
=
\pi(a\xi^{(L_V, M, \pi)}_a(a\bsLV b)(b\bsLV c))
=
\mu(\pi(a), \pi(b), \pi(c)).
\]
Hence, the map $\pi\id_{L_V}^{-1}: (L_V, \mu_{L_V})\to (M, \mu)$
is a homomorphism
\eqref{eqn:characterization:systemhom}
of ternary systems,
and consequently,
$(L_V, (L_V, \mu_{L_V}), \id_{L_V})\sim(L_V, M, \pi)$.
\end{proof}
\section{Examples of ternary systems}
\label{sect:examples}
This section describes several ternary systems
(Definition \ref{defn:construction:ternarysys})
satisfying
\eqref{eqn:construction:defM1}
and
\eqref{eqn:construction:defM2}.
Later we shall characterize
the dynamical YB maps $\Rlmpi(\lambda)$
\eqref{eqn:construction:Rlmpi} constructed by means of
these ternary systems.
\begin{exmp}\label{exmp:examples:fafc}
Let $M$ be a nonempty set, and
$f$ a map from the set $M$ to $M$.
We define the ternary operations on $M$ by:
\begin{eqnarray*}
&&
\mu(a, b, c)=f(a)\quad(\forall a, b, c\in M);
\\
&&
\mu(a, b, c)=f(c)\quad(\forall a, b, c\in M).
\end{eqnarray*}
Each ternary system $(M, \mu)$ defined above satisfies
\eqref{eqn:construction:defM1}
and
\eqref{eqn:construction:defM2}.
\end{exmp}
\begin{rem}
Example \ref{exmp:examples:fafc} satisfying $f=\id_M$ produces degenerate
YB maps in \cite{bukhshtaber98}.
\begin{enumerate}
\item
Let $\mu$ denote the ternary operation on $L$ defined by
$\mu(a, b, c)=c$.
If $L$ is a left quasigroup together with the binary operation $uv:=v$,
then $R^{(L, (L, \mu), \id_L)}(\lambda)(u, v)=(v, v)$
$(\lambda, u, v\in L)$.
This is the map $\Perm{L}\Delta_2$ in \cite{bukhshtaber98}.
Here $\Perm{L}$ is the map \eqref{eqn:construction:flip}.
\item
If $L=(L, \cdot, e_L)$ is a group
and
$\mu(a, b, c)=c$ $(a, b, c\in L)$,
then $\Rlmpi(\lambda)$ is the map $\Perm{L}\mu_1$
in \cite{bukhshtaber98}.
If $L=(L, \cdot, e_L)$ is an abelian group
and
$\mu(a, b, c)=a$ $(a, b, c\in L)$,
then $\Rlmpi(\lambda)$ is the map $\Perm{L}\mu_2$
in \cite{bukhshtaber98}.
\end{enumerate}
\end{rem}
\begin{exmp}\label{exmp:examples:fb}
Let $M$ be a nonempty set, and
$f$ a map from $M$ to $M$
satisfying $f^2=f$.
We define the ternary operation on $M$ by
\[
\mu(a, b, c)=f(b)\quad(\forall a, b, c\in M).
\]
This ternary system $(M, \mu)$ satisfies
\eqref{eqn:construction:defM1}
and
\eqref{eqn:construction:defM2}.
\end{exmp}
\begin{exmp}
Let $M_1=(M_1, \mu_1)$ and $M_2=(M_2, \mu_2)$ be ternary systems
satisfying 
\eqref{eqn:construction:defM1}
and
\eqref{eqn:construction:defM2}.
We denote by $M$ the direct product $M_1\times M_2$
of the sets $M_1$ and $M_2$;
in addition,
let us define the ternary operation $\mu$ on the set $M$ by
\begin{eqnarray*}
&&\mu(a, b, c)
=
(\mu_1(a_1, b_1, c_1),
\mu_2(a_2, b_2, c_2))
\\
&&(a=(a_1, a_2), b=(b_1, b_2), c=(c_1, c_2)\in M=M_1\times M_2).
\end{eqnarray*}
This ternary system $(M, \mu)$ satisfies
\eqref{eqn:construction:defM1}
and
\eqref{eqn:construction:defM2}.
\end{exmp}

We shall introduce three ternary operations
$\mu^G_1$ \eqref{eqn:examples:muG1},
$\mu^G_2$ \eqref{eqn:examples:muG2},
and $\mu^G_3$ \eqref{eqn:examples:muG3}
produced by left quasigroups.

Let $G=(G, *)$ be a left quasigroup
(see Definition \ref{defn:construction:leftquasigroup})
satisfying
that
\begin{equation}\label{eqn:examples:leftquasi1}
(a*c)\bsG((a*b)*c)=(a'*c)\bsG((a'*b)*c)
\quad(\forall a, a', b, c\in G).
\end{equation}
Groups, the quasigroup $(\{ 1, 2, 3\}, *)$ in Example
\ref{exmp:construction:123},
and
the left quasigroups having the right distributive law
\[
(x*y)*z=(x*z)*(y*z)
\quad(\forall x, y, z\in G)
\]
satisfy \eqref{eqn:examples:leftquasi1}
(see below the proof of Proposition \ref{prop:unitary:example123}).
For distributive quasigroups, see \cite[Section V.2]{pflugfelder}.

We define the ternary operations $\mu^G_1$
and
$\mu^G_2$
on the left quasigroup $G$
by:
\begin{eqnarray}
&&\label{eqn:examples:muG1}
\mu^G_1(a, b, c)=a*(b\bsG c)
\quad(a, b, c\in G);
\\\label{eqn:examples:muG2}
&&
\mu^G_2(a, b, c)=c*(b\bsG a)
\quad(a, b, c\in G).
\end{eqnarray}
Here $\bsG$ is the left division \eqref{eqn:construction:leftdivision}
on $G$.
\begin{prop}\label{prop:examples:ternaryG12}
These ternary systems $(G, \mu^G_1)$
and
$(G, \mu^G_2)$
satisfy
$\eqref{eqn:construction:defM1}$
and
$\eqref{eqn:construction:defM2}$.
\end{prop}
\begin{proof}
We shall only prove that the ternary system
$(G, \mu^G_1)$
satisfies
\eqref{eqn:construction:defM1}
and
\eqref{eqn:construction:defM2}.

The following lemma gives rise to \eqref{eqn:construction:defM1}.
Its proof is immediate from \eqref{eqn:examples:muG1}.
\begin{lem}\label{lem:examples:obcA11}
For $a, b, c, d\in G$,
$\mu^G_1(a, b, \mu^G_1(b, c, d))=\mu^G_1(a, c, d)$.
\end{lem}

We shall prove \eqref{eqn:construction:defM2}.
In view of \eqref{eqn:examples:muG1},
\begin{eqnarray*}
&&
\mbox{LHS of \eqref{eqn:construction:defM2}}
\\
&=&
(a*(b\bsG c))*(c\bsG d)
\\
&=&
(a*(c\bsG d))*((a*(c\bsG d))\bsG ((a*(b\bsG c))*(c\bsG d))),
\\
&&
\mbox{RHS of \eqref{eqn:construction:defM2}}
\\
&=&
(a*(c\bsG d))*((b*(c\bsG d))\bsG d)
\\
&=&
(a*(c\bsG d))*((b*(c\bsG d))\bsG ((b*(b\bsG c))*(c\bsG d))).
\end{eqnarray*}
The right-hand-sides of the above equations are the same,
because of \eqref{eqn:examples:leftquasi1}.
\end{proof}
\begin{rem}\label{rem:examples:shibukawa05}
Every dynamical YB map $R^{(G)}(\lambda)$
\eqref{eqn:summary:RLP}
constructed in the work
\cite{shibukawa05}
is produced by the ternary system
$(G, \mu^G_1)$ \eqref{eqn:examples:muG1}.
Let $L$ be a loop,
$G=(G, *, e_G)$ a group, and $\pi:L\to G$ a bijection
satisfying $\pi(e_L)=e_G$.
Here $e_G$ is the unit element of the group $G$.
By the definitions \eqref{eqn:summary:theta}
of the maps $\theta(u)$ and $\theta(u)^{-1}$,
the map $\xi_\lambda^{(G)}(u)$ \eqref{eqn:summary:xishibukawa05}
($\lambda, u\in L$)
is expressed as
\begin{equation}\label{eqn:examples:shibukawa05}
\xi_\lambda^{(G)}(u)(v)
=
\lambda\bsL\pi^{-1}(\mu^G_1(\pi(\lambda), \pi(\lambda u), \pi((\lambda u)v)))
\quad(v\in L).
\end{equation}
On account of \eqref{eqn:summary:xishibukawa05},
\eqref{eqn:summary:etashibukawa05},
\eqref{eqn:construction:xilmpi},
and
\eqref{eqn:construction:etalmpi},
\eqref{eqn:examples:shibukawa05}
induces that
all the dynamical YB maps $R^{(G)}(\lambda)$
are constructed by means of the ternary systems $(G, \mu^G_1)$.
\end{rem}

Next task is to define the ternary operation $\mu^G_3$
\eqref{eqn:examples:muG3}.

Let $G=(G, *)$ be a left quasigroup.
We suppose that $G$ satisfies the following for all $a, b, c, d\in G$:
{\allowdisplaybreaks
\begin{eqnarray}
\label{eqn:examples:leftquasi22}
&&
(b*c)*(a\bsG ((a*c)*((b*c)\bsG (b*d))))
\\\nonumber
&=&
b*(a\bsG ((a*c)*d));
\\
&&\label{eqn:examples:leftquasi21}
(a*c)*((b*c)\bsG (b*d))
\\\nonumber
&=&
((a*c)*d)*((b*(a\bsG ((a*c)*d)))\bsG (b*d)).
\end{eqnarray}
}
If $G$ is a group,
$G$ satisfies
\eqref{eqn:examples:leftquasi22}
and
\eqref{eqn:examples:leftquasi21}.

Let $\mu^G_3$ denote the ternary operation
on the set $G$ defined by
\begin{equation}\label{eqn:examples:muG3}
\mu^G_3(a, b, c)=b*(a\bsG c)
\quad(a, b, c\in G).
\end{equation}
\begin{prop}\label{prop:examples:ternaryG3}
This ternary system $(G, \mu^G_3)$
satisfies
$\eqref{eqn:construction:defM1}$
and
$\eqref{eqn:construction:defM2}$.
\end{prop}
\begin{proof}
We shall only prove that the ternary system $(G, \mu^G_3)$
satisfies \eqref{eqn:construction:defM1}.
In view of \eqref{eqn:examples:muG3},
\begin{eqnarray*}
&&
\mbox{LHS of \eqref{eqn:construction:defM1}}
\\
&=&
(b*(a\bsG c))*(a\bsG (c*((b*(a\bsG c))\bsG d)));
\\
&&
\mbox{RHS of \eqref{eqn:construction:defM1}}
\\
&=&
b*(a\bsG (c*(b\bsG d)))
\\
&=&
b*(a\bsG ((a*(a\bsG c))*(b\bsG d))).
\end{eqnarray*}
The right-hand-sides of the above equations are the same,
because of \eqref{eqn:examples:leftquasi22}.
\end{proof}

Final task in this section is to
characterize the dynamical YB map
$\Rlmpi(\lambda)$ \eqref{eqn:construction:Rlmpi}
that
the ternary systems \eqref{eqn:examples:muG1},
\eqref{eqn:examples:muG2},
and \eqref{eqn:examples:muG3}
define.

Let $\cA_1$ denote the subcategory of the category
$\cA$ whose objects and morphisms are defined as follows:
$V\in Ob(\cA)$ is an object of $\cA_1$,
iff there exists a representative
$(L, (M, \mu), \pi)$ of $V$ such that
the ternary operation $\mu$ on $M$ satisfies
\begin{eqnarray}
&&\label{eqn:examples:obA11}
\mu(a, b, \mu(b, c, d))=\mu(a, c, d)
\quad(\forall a, b, c, d\in M),
\\\label{eqn:examples:obA12}
&&
\mu(a, a, b)=b
\quad(\forall a, b\in M);
\end{eqnarray}
$f: V\to V'\in \Hom(\cA)$ is a morphism of $\cA_1$,
iff $V, V'\in Ob(\cA_1)$.

Let $L$ be a left quasigroup,
$G$ a left quasigroup satisfying $\eqref{eqn:examples:leftquasi1}$,
and $\pi$ a $($set-theoretical\/$)$ bijection from $L$ to $G$.
$[(L, (G, \mu^G_1), \pi)]$ is an object of the category $\cA$
because
of Proposition
\ref{prop:examples:ternaryG12}.
Moreover,
\begin{prop}\label{prop:examples:LGpicA1}
$[(L, (G, \mu^G_1), \pi)]$ is an object of the category $\cA_1$.
\end{prop}
\begin{proof}
The proof is immediate from
\eqref{eqn:examples:muG1}
(see Lemma \ref{lem:examples:obcA11}).  
\end{proof}

Conversely, every object of the category $\cA_1$ is expressed by means of
the
ternary system $(G, \mu^G_1)$.
\begin{prop}\label{prop:examples:12}
If $V\in Ob(\cA_1)$,
then
there exist a left quasigroup
$(G, *)$ satisfying $\eqref{eqn:examples:leftquasi1}$
and a bijection $\pi': L_V\to G$
such that
$V=[(L_V, (G, \mu^G_1), \pi')]$.
\end{prop}
\begin{proof}
We denote by $L$ the left quasigroup $L_V$.
Since $V\in Ob(\cA_1)$,
there exists a representative
$(L, (M, \mu), \pi)$ of $V$ such that
the ternary operation $\mu$ on $M$ satisfies
\eqref{eqn:examples:obA11}
and
\eqref{eqn:examples:obA12}.

We fix any element $\lambda\in L$.
Let $\alpha$ and $\beta$ denote
the following binary operations on $L$:
for $u, v\in L$,
\begin{eqnarray}
&&\label{eqn:examples:alpha}
\alpha(u, v)=\lambda\bsL\pi^{-1}(\mu(\pi(\lambda), \pi(\lambda u),
\pi(\lambda v)));
\\\label{eqn:examples:beta}
&&
\beta(u, v)=\lambda\bsL\pi^{-1}(\mu(\pi(\lambda u), \pi(\lambda),
\pi(\lambda v))).
\end{eqnarray}
\begin{lem}
For all $u, v\in L$,
$\beta(u, \alpha(u, v))=v$
and
$\alpha(u, \beta(u, v))=v$.
\end{lem}
\begin{proof}
The proof is immediate from \eqref{eqn:examples:obA11}
and
\eqref{eqn:examples:obA12}.
\end{proof}

We denote by $G$ and $(*)$ the set $L$ and the 
binary operation $\beta$ on the set $G(=L)$, respectively.
The above lemma gives rise to that
$G=(G, *)$ is a left quasigroup.
The left division
on $G$ is the binary operation $\alpha$;
$a\bsG c=\alpha(a, c)$ $(a, c\in G)$.
\begin{lem}\label{lem:examples:leftquasi1}
The left quasigroup $G$ satisfies
$\eqref{eqn:examples:leftquasi1}$.
\end{lem}
\begin{proof}
Because of \eqref{eqn:examples:obA11}, \eqref{eqn:examples:alpha},
and \eqref{eqn:examples:beta},
\begin{equation}\label{eqn:examples:abc}
a*(b\bsG c)=
\lambda\bsL\pi^{-1}(\mu(\pi(\lambda a), \pi(\lambda b), \pi(\lambda c)))
\quad(\forall a, b, c\in G).
\end{equation}

Let $a$, $a'$, $b$, and $c$ be elements of $G(=L)$.
In view of \eqref{eqn:examples:abc},
{\allowdisplaybreaks
\begin{eqnarray*}
&&(a*b)*c
\\
&=&
(a*(a'\bsG (a'*b)))*((a'*b)\bsG((a'*b)*c))
\\
&=&
\lambda\bsL\pi^{-1}(\mu(\mu(\pi(\lambda a), \pi(\lambda a'),
\pi(\lambda (a'*b))), \pi(\lambda (a'*b)), \pi(\lambda ((a'*b)*c))));
\\
&&
(a*c)*((a'*c)\bsG ((a'*b)*c))
\\
&=&
(a*(a'\bsG(a'*((a'*b)\bsG ((a'*b)*c)))))*
\\
&&
*((a'*((a'*b)\bsG ((a'*b)*c)))\bsG ((a'*b)*c))
\\
&=&
\lambda\bsL\pi^{-1}(\mu(
\mu(\pi(\lambda a), \pi(\lambda a'),
\mu(\pi(\lambda a'), \pi(\lambda (a'*b)), \pi(\lambda ((a'*b)*c)))),
\\
&&
\mu(\pi(\lambda a'), \pi(\lambda (a'*b)), \pi(\lambda ((a'*b)*c))),
\pi(\lambda ((a'*b)*c)))).
\end{eqnarray*}
}

With the aid of \eqref{eqn:construction:defM2},
\[(a*b)*c
=
(a*c)*((a'*c)\bsG ((a'*b)*c))
\]
for all $a, a', b, c\in G$.
This is equivalent to \eqref{eqn:examples:leftquasi1}.
\end{proof}
Let $\pi'$ denote the map from $L$ to $G(=L)$ defined by
$\pi'(u)=\lambda\bsL u$ $(u\in L)$.
\begin{lem}
The map $\pi'$ is bijective\/$;$
$\pi'^{-1}(a)=\lambda a$ $(a\in G)$.
\end{lem}

Finally, we shall demonstrate that $V=[(L, (G, \mu^G_1), \pi')]$.
From \eqref{eqn:examples:abc},
\[
\pi'\pi^{-1}(\mu(\pi\pi'^{-1}(a), \pi\pi'^{-1}(b), \pi\pi'^{-1}(c)))
=a*(b\bsG c)
=\mu^G_1(a, b, c).
\]
Hence, the map $\pi\pi'^{-1}: (G, \mu^G_1)\to M$ is a
homomorphism \eqref{eqn:characterization:systemhom}
of ternary systems.
As a result,
$(L, M, \pi)\sim (L, (G, \mu^G_1), \pi')$;
that is, $V=[(L, (G, \mu^G_1), \pi')]$.

This completes the proof of Proposition \ref{prop:examples:12}.
\end{proof}

We denote by $\cD_1$ the subcategory of the category $\cD$
whose objects and morphisms are defined as follows:
$V=(L, R(\lambda))\in Ob(\cD)$ is an object of $\cD_1$,
iff the map $\xi_\lambda(u)$
$(\lambda, u\in L)$
\eqref{eqn:characterization:etaxi} satisfies
\begin{eqnarray*}
&&
\xi_\lambda(u)\xi_{\lambda u}(v)=\xi_\lambda(\lambda\bsL ((\lambda u)v))
\quad(\forall \lambda, u, v\in L),
\\
&&
\xi_\lambda(\lambda\bsL\lambda)=\id_L
\quad(\forall \lambda\in L);
\end{eqnarray*}
$f: V\to V'\in \Hom(\cD)$ is a morphism of $\cD_1$,
iff $V, V'\in Ob(\cD_1)$.

The functors $S: \cA\to \cD$ and $T: \cD\to \cA$ induce the following. 
\begin{prop}\label{prop:examples:isomorphic1}
The category $\cA_1$ is isomorphic to the category $\cD_1$.
\end{prop}
\begin{proof}
The proof is straightforward.
\end{proof}

Let us introduce subcategories $\cA_2$
and $\cA_3$ (resp.\ $\cD_2$ and $\cD_3$)
of the category
$\cA$ (resp.\ $\cD$),
which characterize the dynamical YB maps
$\Rlmpi(\lambda)$
\eqref{eqn:construction:Rlmpi}
constructed by means of the ternary systems
\eqref{eqn:examples:muG2} and
\eqref{eqn:examples:muG3}:
$V\in Ob(\cA)$ is an object of $\cA_2$,
iff there exists a representative
$(L, (M, \mu), \pi)$ of $V$ such that
the ternary operation $\mu$ on $M$ satisfies
\begin{eqnarray}\label{eqn:examples:obcA21}
&&
\mu(\mu(a, b, c), c, d)=\mu(a, b, d)
\quad(\forall a, b, c, d\in M),
\\\label{eqn:examples:obcA22}
&&
\mu(a, b, b)=a
\quad(\forall a, b\in M);
\end{eqnarray}
$f: V\to V'\in \Hom(\cA)$ is a morphism of $\cA_2$,
iff $V, V'\in Ob(\cA_2)$;
$V\in Ob(\cA)$ is an object of $\cA_3$,
iff there exists a representative
$(L, (M, \mu), \pi)$ of $V$ such that
the ternary operation $\mu$ on $M$ satisfies
\begin{eqnarray}\label{eqn:examples:obcA31}
&&
\mu(a, b, c)=\mu(d, b, \mu(a, d, c))
\quad(\forall a, b, c, d\in M),
\\\label{eqn:examples:obcA32}
&&
\mu(a, a, b)=b
\quad(\forall a, b\in M);
\end{eqnarray}
$f: V\to V'\in \Hom(\cA)$ is a morphism of $\cA_3$,
iff $V, V'\in Ob(\cA_3)$;
$V=(L, R(\lambda))\in Ob(\cD)$ is an object of $\cD_2$,
iff the maps $\xi_\lambda(u)$ and $\eta_\lambda(v)$
$(\lambda, u, v\in L)$
\eqref{eqn:characterization:etaxi}
satisfy
\begin{eqnarray*}
&&
(\lambda\xi_\lambda(u)(v))
\xi_{\lambda \xi_\lambda(u)(v)}(\eta_\lambda(v)(u))(w)
\\
&=&
\lambda\xi_\lambda(u)((\lambda u)\bsL (((\lambda u)v)w))
\quad(\forall \lambda, u, v, w\in L),
\\
&&
\xi_\lambda(u)((\lambda u)\bsL(\lambda u))=\lambda\bsL\lambda
\quad(\forall \lambda, u\in L);
\end{eqnarray*}
$f: V\to V'\in \Hom(\cD)$ is a morphism of $\cD_2$,
iff $V, V'\in Ob(\cD_2)$;
$V=(L, R(\lambda))\in Ob(\cD)$ is an object of $\cD_3$,
iff the map $\xi_\lambda(u)$
$(\lambda, u\in L)$
\eqref{eqn:characterization:etaxi} satisfies
\begin{eqnarray*}
&&
\lambda\xi_\lambda(v)
((\lambda v)\bsL
((\lambda u)\xi_{\lambda u}
((\lambda u)\bsL\lambda)(w)))
\\
&=&
(\lambda u)\xi_{\lambda u}((\lambda u)\bsL (\lambda v))
((\lambda v)\bsL(\lambda w))
\quad(\forall \lambda, u, v, w\in L),
\\
&&
\xi_\lambda(\lambda\bsL\lambda)=\id_L
\quad(\forall \lambda\in L);
\end{eqnarray*}
$f: V\to V'\in \Hom(\cD)$ is a morphism of $\cD_3$,
iff $V, V'\in Ob(\cD_3)$.

The functors $S: \cA\to \cD$ and $T: \cD\to \cA$ give rise to the following proposition.
\begin{prop}
The categories $\cA_2$ and $\cA_3$
are isomorphic to the categories
$\cD_2$ and $\cD_3$, respectively.
\end{prop}
The proof of the following proposition is immediate from
\eqref{eqn:examples:muG2}.
\begin{prop}
Let $L$ be a left quasigroup,
$G$ a left quasigroup satisfying $\eqref{eqn:examples:leftquasi1}$,
and $\pi$ a $($set-theoretical\/$)$ bijection from $L$ to $G$.
Then $[(L, (G, \mu^G_2), \pi)]$ is an object of the category $\cA_2$.
\end{prop}
\begin{prop}
If $V\in Ob(\cA_2)$,
then
there exist a left quasigroup
$(G, *)$ satisfying $\eqref{eqn:examples:leftquasi1}$
and a bijection $\pi': L_V\to G$
such that
$V=[(L_V, (G, \mu^G_2), \pi')]$.
\end{prop}
\begin{proof}
The proof is similar to that of Proposition
\ref{prop:examples:12}.
For the reason that $V\in Ob(\cA_2)$,
there exists a representative
$(L_V, (M, \mu), \pi)$ of $V$ such that
the ternary operation $\mu$ on $M$ satisfies
\eqref{eqn:examples:obcA21}
and
\eqref{eqn:examples:obcA22}.

Let $G$ denote the set $L_V$.
We fix any element $\lambda\in G(=L_V)$.
Let us define the binary operation $*$ on $G(=L_V)$ by
\[
a*b=\lambda\bsLV\pi^{-1}(\mu(\pi(\lambda b), \pi(\lambda), \pi(\lambda a)))
\quad(a, b\in G).
\]
Then $(G, *)$ is a left quasigroup;
the left division $\bsG$ is as follows.
\[
a\bsG c=\lambda\bsLV\pi^{-1}
(\mu(\pi(\lambda c), \pi(\lambda a), \pi(\lambda)))
\quad(a, c\in G(=L_V)).
\]
The following equation
induces that $(G, *)$ satisfies
\eqref{eqn:examples:leftquasi1}
(cf.\ Lemma \ref{lem:examples:leftquasi1}).
\[
c*(b\bsG a)=
\lambda\bsLV\pi^{-1}(\mu(\pi(\lambda a), \pi(\lambda b), \pi(\lambda c)))
\quad(\forall a, b, c\in G).
\]
Moreover, we denote by $\pi'$ the bijection from $L_V$ to $G(=L_V)$
defined by
$\pi'(u)=\lambda\bsLV u$ $(u\in L_V)$.
The triplet $(L_V, (G, \mu^G_2), \pi')$ is a representative of $V$;
that is, $V=[(L_V, (G, \mu^G_2), \pi')]$.
\end{proof}

Let $L$ be a left quasigroup,
$G$ a left quasigroup satisfying $\eqref{eqn:examples:leftquasi22}$
and
$\eqref{eqn:examples:leftquasi21}$,
and $\pi$ a $($set-theoretical\/$)$ bijection from $L$ to $G$.
Then $[(L, (G, \mu^G_3), \pi)]$ is an object of the category $\cA$
by virtue of Proposition
\ref{prop:examples:ternaryG3}.
\begin{prop}
$[(L, (G, \mu^G_3), \pi)]$ is an object of the category $\cA_3$.
\end{prop}
\begin{prop}
If $V\in Ob(\cA_3)$,
then
there exist a left quasigroup
$(G, *)$ satisfying $\eqref{eqn:examples:leftquasi22}$
and
$\eqref{eqn:examples:leftquasi21}$,
and a bijection $\pi': L_V\to G$
such that
$V=[(L_V, (G, \mu^G_3), \pi')]$.
\end{prop}
\begin{proof}
Because $V\in Ob(\cA_3)$,
there exists a representative
$(L_V, (M, \mu), \pi)$ of $V$ such that
the ternary operation $\mu$ on $M$ satisfies
\eqref{eqn:examples:obcA31}
and
\eqref{eqn:examples:obcA32}.

Let $G$ denote the set $L_V$.
We fix any element $\lambda\in G(=L_V)$.
Let us define the binary operation $*$ on $G(=L_V)$ by
\[
a*b=\lambda\bsLV\pi^{-1}(\mu(\pi(\lambda), \pi(\lambda a), \pi(\lambda b)))
\quad(a, b\in G).
\]
Then $(G, *)$ is a left quasigroup satisfying
\eqref{eqn:examples:leftquasi22}
and
\eqref{eqn:examples:leftquasi21};
the left division $\bsG$ is as follows.
\[
a\bsG c=\lambda\bsLV\pi^{-1}
(\mu(\pi(\lambda a), \pi(\lambda), \pi(\lambda c)))
\quad(a, c\in G(=L_V)).
\]
The proof of \eqref{eqn:examples:leftquasi22}
and
\eqref{eqn:examples:leftquasi21}
is due to
the following
(cf.\ Lemma \ref{lem:examples:leftquasi1}).
\[
b*(a\bsG c)=
\lambda\bsLV\pi^{-1}(\mu(\pi(\lambda a), \pi(\lambda b), \pi(\lambda c)))
\quad(\forall a, b, c\in G).
\]
In addition, we denote by $\pi'$ the bijection from $L_V$ to $G(=L_V)$
defined by
$\pi'(u)=\lambda\bsLV u$ $(u\in L_V)$.

These
$(G, *)$ and $\pi'$ are what we desire.
\end{proof}
\section{Unitary condition}
\label{sect:unitary}
Let $(L, M, \pi)$ be an element of $LMB$.
In this section, we discuss the unitary condition
\eqref{eqn:summary:unitary}
of the dynamical YB map $\Rlmpi(\lambda)$
\eqref{eqn:construction:Rlmpi}.
\begin{prop}\label{prop:unitary:unitary}
The dynamical YB map $\Rlmpi(\lambda)$
satisfies the unitary condition,
if and only if
the ternary operation $\mu$ on $M$ satisfies that
\begin{equation}\label{eqn:unitary:unitary}
\mu(a, \mu(a, b, c), c)=b
\quad(\forall a, b, c\in M).
\end{equation}
\end{prop}
\begin{proof}
Let $\lambda$ be an element of the left quasigroup $L$.
We define the map $\tilde{f}_{\lambda}: L\times L\to L\times L$ by
\begin{equation}
\label{eqn:unitary:tildeflambda}
\tilde{f}_{\lambda}(u, v)
=(\lambda u, (\lambda u) v)
\quad(u, v\in L).
\end{equation}
\begin{lem}\label{lem:unitary:f}
The map $\tilde{f}_\lambda$ is bijective\/$;$
$\tilde{f}^{-1}_\lambda(u, v)=(\lambda\bsL u,
u\bsL v)$
$(u, v\in L)$.
\end{lem}
By using the maps
$\Perm{L}$
\eqref{eqn:construction:flip},
$s(\pi(\lambda))$
\eqref{eqn:construction:defsa},
$\tilde{f}_\lambda$
\eqref{eqn:unitary:tildeflambda},
and
$\tilde{f}^{-1}_\lambda$,
the dynamical YB map $\Rlmpi(\lambda)$ is expressed
as follows.
\begin{equation}
\label{eqn:unitary:Rlmpi}
\Rlmpi(\lambda)=\Perm{L}\tilde{f}^{-1}_\lambda(\pi^{-1}\times \pi^{-1})
s(\pi(\lambda))(\pi\times\pi)\tilde{f}_\lambda.
\end{equation}
Since the
map $\pi: L\to M$
is bijective,
the unitary condition of the dynamical YB map $\Rlmpi(\lambda)$ is equivalent
to that $s(a)^2=\id_{M\times M}$ for all $a\in M$.

The rest of the proof is immediate from the definition
\eqref{eqn:construction:defsa}
of the map $s(a)$.
\end{proof}
Let $M$ denote one of the ternary systems in Section
\ref{sect:examples}.
By using Proposition
\ref{prop:unitary:unitary},
we shall clarify a necessary and sufficient condition for the dynamical YB map $\Rlmpi(\lambda)$
\eqref{eqn:construction:Rlmpi}
to satisfy the unitary condition.

Let us suppose that $M$ is a ternary system in Example
\ref{exmp:examples:fafc}.
From Proposition \ref{prop:unitary:unitary},
the dynamical YB map $\Rlmpi(\lambda)$
does not satisfy the unitary condition, unless $|M|=1$.
If $|M|=1$,
then the dynamical YB map $\Rlmpi(\lambda)=\id_{L\times L}$;
hence, it is trivial that the dynamical YB map $\Rlmpi(\lambda)$
satisfies the unitary condition.

Next we suppose that $M$ is a ternary system in Example
\ref{exmp:examples:fb}.
The dynamical YB map $\Rlmpi(\lambda)$
satisfies the unitary condition, if and only if
the map $f$ is the identity map $\id_M$.
If $f=\id_M$, then $\Rlmpi(\lambda)=\Perm{L}$
\eqref{eqn:construction:flip}.
\begin{prop}\label{prop:unitary:example123}
Let $M$ be a ternary system 
$(G, \mu^G_1)$ $\eqref{eqn:examples:muG1}$
or $(G, \mu^G_2)$ $\eqref{eqn:examples:muG2}$.
The dynamical YB map $\Rlmpi(\lambda)$
satisfies the unitary condition, if and only if
\begin{equation}\label{eqn:unitary:exmp12}
(a*b)*c=(a*c)*b
\quad(\forall a, b, c\in G).
\end{equation}
\end{prop}
\begin{proof}
Let $G=(G, *)$ be a left quasigroup satisfying
\eqref{eqn:examples:leftquasi1}.
We shall only show this proposition in the case that $M=(G, \mu^G_1)$.

Let us suppose that the dynamical YB map $\Rlmpi(\lambda)$ satisfies the unitary
condition.
It follows from
\eqref{eqn:examples:muG1}
and
Proposition \ref{prop:unitary:unitary}
that
\begin{eqnarray*}
a*((a*c)\bsG ((a*b)*c))&=&\mu^G_1(a, \mu^G_1(a, a*b, (a*b)*c), (a*b)*c)
\\
&=&a*b
\qquad(\forall a, b, c\in G).
\end{eqnarray*}
This is equivalent to \eqref{eqn:unitary:exmp12},
since $G$ is a left quasigroup.

Conversely, we suppose that Eq.\ \eqref{eqn:unitary:exmp12}
holds.
With the aid of \eqref{eqn:examples:muG1},
\begin{eqnarray}\label{eqn:unitary:proof}
\mbox{LHS of \eqref{eqn:unitary:unitary}}
&=&
a*((a*(b\bsG c))\bsG c)
\\\nonumber
&=&
a*((a*(b\bsG c))\bsG
((a*(a\bsG b))*(b\bsG c)))
\end{eqnarray}
for all $a, b, c\in G$.
By virtue of \eqref{eqn:unitary:exmp12},
\begin{eqnarray*}
\mbox{RHS of \eqref{eqn:unitary:proof}}
&=&
a*((a*(b\bsG c))\bsG
((a*(b\bsG c))*(a\bsG b)))
\\
&=&
b.
\end{eqnarray*}
Because of Proposition \ref{prop:unitary:unitary},
the dynamical YB map $\Rlmpi(\lambda)$
satisfies the unitary condition.
\end{proof}
Let $G$ be an abelian group or the quasigroup $(\{ 1, 2, 3\}, *)$
in Example \ref{exmp:construction:123}.
This $G$ satisfies \eqref{eqn:unitary:exmp12}
(the proof is straightforward);
in addition,
\eqref{eqn:unitary:exmp12} induces that
$G$ satisfies
\eqref{eqn:examples:leftquasi1}.
\begin{rem}
If the left quasigroup $G$ is a group, then
\eqref{eqn:unitary:exmp12}
is equivalent to that
the group $G$ is abelian.
Hence, Remark \ref{rem:examples:shibukawa05}
and Proposition \ref{prop:unitary:example123}
reproduce Theorem
\ref{thm:summary:unitary}.
\end{rem}

The proof of the proposition below is similar to that of Proposition
\ref{prop:unitary:example123}.
\begin{prop}\label{prop:unitary:example123-2}
Let $M$ be a ternary system $(G, \mu^G_3)$
$\eqref{eqn:examples:muG3}$.
The dynamical YB map $\Rlmpi(\lambda)$
satisfies the unitary condition, if and only if
\[
(a*b)*b=a
\quad(\forall a, b\in G).
\]
\end{prop}
\section{IRF-IRF correspondence}
\label{sect:correspondence}
Let $L_i=(L_i, \bullet_i)$
$(i=1, 2)$ be left quasigroups
(see Definition \ref{defn:construction:leftquasigroup}),
$M=(M, \mu)$ a ternary system
(Definition \ref{defn:construction:ternarysys})
satisfying \eqref{eqn:construction:defM1}
and
\eqref{eqn:construction:defM2},
and
$\pi_i: L_i\to M$
$(i=1, 2)$ bijections.
Let $\lambda_1$ and $\lambda_2$ be elements of the left quasigroups $L_1$
and $L_2$, respectively.

To end this paper, we establish a correspondence between
two dynamical YB maps
$R^{(L_1, M, \pi_1)}(\lambda_1)$
and
$R^{(L_2, M, \pi_2)}(\lambda_2)$
\eqref{eqn:construction:Rlmpi}
called an {\em IRF-IRF correspondence}
(Proposition \ref{prop:correspondence:IRFIRF}).
Let $\tilde{f}^{(i)}_{\lambda_i}$ $(i=1, 2)$ denote the map
from $L_i\times L_i$ to itself defined by
\eqref{eqn:unitary:tildeflambda}.
By means of these maps, we define the map
$J(\lambda_1): L_1\times L_1\to L_2\times L_2$
as follows.
\[
J(\lambda_1)
=
\tilde{f}^{(2)\,-1}_{\pi^{-1}_2\pi_1(\lambda_1)}
(\pi^{-1}_2\pi_1\times\pi^{-1}_2\pi_1)
\tilde{f}^{(1)}_{\lambda_1}\Perm{L_1}.
\]
Here
$\Perm{L_1}$ is the map
\eqref{eqn:construction:flip}.
Let us define the maps
$R^{(L_2, M, \pi_2)\, 21}(\lambda_2): L_2\times L_2\to L_2\times L_2$
and
$J^{21}(\lambda_1): L_1\times L_1\to L_2\times L_2$
by
$R^{(L_2, M, \pi_2)\, 21}(\lambda_2)=\Perm{L_2}
R^{(L_2, M, \pi_2)}(\lambda_2)\Perm{L_2}$
and
$J^{21}(\lambda_1)
=
\Perm{L_2}
J(\lambda_1)
\Perm{L_1}$.

Eq.\ 
\eqref{eqn:unitary:Rlmpi}
plays an essential role in the proof of Proposition
\ref{prop:correspondence:IRFIRF}.
\begin{prop}\label{prop:correspondence:IRFIRF}
$R^{(L_1, M, \pi_1)}(\lambda_1)
=
J(\lambda_1)^{-1}
R^{(L_2, M, \pi_2)\, 21}(\pi^{-1}_2\pi_1(\lambda_1))
J^{21}(\lambda_1).
$
\end{prop}
\begin{proof}
Let $a$ denote the element $\pi_1(\lambda_1)$
of the ternary system $M$.
From \eqref{eqn:unitary:Rlmpi},
\[
s(a)
=
(\pi_i\times \pi_i)
\tilde{f}^{(i)}_{\pi^{-1}_i(a)}\Perm{L_i}
R^{(L_i, M, \pi_i)}(\pi^{-1}_i(a))
\tilde{f}^{(i)\,-1}_{\pi^{-1}_i(a)}
(\pi^{-1}_i\times\pi^{-1}_i)
\]
for $i=1, 2$.
This equation immediately
induces Proposition \ref{prop:correspondence:IRFIRF}.
\end{proof}

This IRF-IRF correspondence is said to be a {\em vertex-IRF correspondence},
iff the dynamical YB map $R^{(L_2, M, \pi_2)}(\lambda_2)$
is independent of the dynamical parameter
$\lambda_2$;
hence, $R^{(L_2, M, \pi_2)}(\lambda_2)$
is a YB map.
We denote by $R^{(L_2, M, \pi_2)}$ the YB map
$R^{(L_2, M, \pi_2)}(\lambda_2)$.
The vertex-IRF correspondence is as follows
(cf.\ 
\cite[Definition 5.4]{etingofbook}
and
\cite[(4.10)]{koornwinder}).
\begin{equation}\label{eqn:correspondence:vertexirf}
R^{(L_1, M, \pi_1)}(\lambda_1)
=
J(\lambda_1)^{-1}
R^{(L_2, M, \pi_2)\, 21}
J^{21}(\lambda_1).
\end{equation}
\begin{rem}\label{rem:correspondence:exchange}
The maps $J(\lambda_1)$ and $R^{(L_1, M, \pi_1)}(\lambda_1)$
correspond to the fusion matrix and the exchange matrix,
respectively
(see
\cite[Sections 5.1 and 5.2]{etingofbook}
and
\cite[Sections 3 and 4]{koornwinder}).
\end{rem}
Let $V=(L, R(\lambda))$ be an object of the category $\cD_1$.
We shall discuss a vertex-IRF correspondence whose IRF part is
this dynamical YB map $R(\lambda)$. 
\begin{prop}\label{prop:correspondence:vertexIRF}
The dynamical YB map $R(\lambda)$ has at least one
vertex-IRF correspondence.
\end{prop}
\begin{proof}
From Propositions
\ref{prop:examples:12}
and
\ref{prop:examples:isomorphic1},
there exist a left quasigroup
$(G, *)$ satisfying $\eqref{eqn:examples:leftquasi1}$
and a bijection $\pi: L\to G$
such that
$R(\lambda)=R^{(L, (G, \mu^G_1), \pi)}(\lambda)$
\eqref{eqn:construction:Rlmpi}.

Let $\circ$ denote the binary operation on the set $L$
defined by
$u\circ v=\pi^{-1}(\pi(u)*\pi(v))$
$(u, v\in L)$.
\begin{lem}
$L'=(L, \circ)$ is a left quasigroup\/$;$
the left division $\bsl{L'}$
$\eqref{eqn:construction:leftdivision}$
of $L'$ is defined by
$u\bsl{L'}w=\pi^{-1}(\pi(u)\bsG\pi(w))$
$(u, w\in L')$.
\end{lem}
\begin{proof}
The proof is straightforward, because $G$ is a left quasigroup
and $\pi$ is a bijection.
\end{proof}
By virtue of \eqref{eqn:construction:xilmpi},
\eqref{eqn:examples:muG1},
and
the definition of the binary operation $\circ$,
\begin{eqnarray}\label{eqn:correspondence:xi}
&&
\xi^{(L', (G, \mu^G_1), \pi)}_\lambda(u)(v)
\\\nonumber
&=&
\lambda\bsl{L'}\pi^{-1}(\mu^G_1(\pi(\lambda), \pi(\lambda\circ u),
\pi((\lambda\circ u)\circ v)))
\\\nonumber
&=&
v
\end{eqnarray}
for all $\lambda, u, v\in L'$.
With the aid of \eqref{eqn:construction:etalmpi}
and
\eqref{eqn:correspondence:xi},
\[
\eta^{(L', (G, \mu^G_1), \pi)}_\lambda(v)(u)
=\pi^{-1}((\pi(\lambda)*\pi(v))\bsG((\pi(\lambda)*\pi(u))*\pi(v)))
\]
for all $\lambda, u, v\in L'$.
This element $\eta^{(L', (G, \mu^G_1), \pi)}_\lambda(v)(u)$
is independent of $\lambda$ because of \eqref{eqn:examples:leftquasi1}.
Hence,
the dynamical YB map
$R^{(L', (G, \mu^G_1), \pi)}(\lambda)$
is a YB map.

From Proposition
\ref{prop:correspondence:IRFIRF},
the dynamical YB map $R(\lambda)=R^{(L, (G, \mu^G_1), \pi)}(\lambda)$
has a vertex-IRF correspondence whose vertex counterpart is
this YB map $R^{(L', (G, \mu^G_1), \pi)}(\lambda)$.
\end{proof}
We can construct objects $V=(L, R(\lambda))$ of the category $\cD_1$
such that each $R(\lambda)$ really depends on the
dynamical parameter $\lambda$. 
Let $G=(G, *)$ be a quasigroup satisfying
\eqref{eqn:examples:leftquasi1}
(see Definition \ref{defn:construction:quasigroup}),
and
let $\circ$ denote the following binary operation on $G$:
$a\circ b=b$ $(a, b\in G)$.
Then $L'=(G, \circ)$ is a left quasigroup,
and \eqref{eqn:construction:Rlmpi} gives rise to that
\begin{equation}\label{eqn:correspondence:dependdYB}
R^{(L', (G, \mu^G_1), \id_G)}(\lambda)(u, v)=(v, \lambda*(u\bsG v))
\quad(\forall\lambda, u, v\in G).
\end{equation}
On account of
Propositions
\ref{prop:examples:LGpicA1}
and \ref{prop:examples:isomorphic1},
the pair $(L', R^{(L', (G, \mu^G_1), \id_G)}(\lambda))$ is an object of
the category $\cD_1$.

We suppose that $R^{(L', (G, \mu^G_1), \id_G)}(\lambda)
=
R^{(L', (G, \mu^G_1), \id_G)}(\lambda')$.
Then $\lambda*(u\bsG v)=\lambda'*(u\bsG v)$ for all $u, v\in G$.
This equation induces that $\lambda=\lambda'$,
because $G$ is a quasigroup.
Hence, the dynamical YB map $R^{(L', (G, \mu^G_1), \id_G)}(\lambda)$
is dependent on $\lambda$,
unless $|G|=1$.
\begin{rem}
Let $L$ and $G$ be groups,
and $\pi$ a bijective 1-cocycle of the group $L$ with coefficients in the
group $G$.
From Remarks \ref{rem:summary:lu}
and \ref{rem:examples:shibukawa05},
$R^{(L, (G, \mu^G_1), \pi)}(\lambda)$
is the YB map
in \cite{lu00}.
Let us suppose that $|L|\neq 1$ (hence $|G|\neq 1$).
The YB map
$R^{(L, (G, \mu^G_1), \pi)}(\lambda)$
has a vertex-IRF correspondence whose IRF
counterpart is the dynamical YB map
$R^{(L', (G, \mu^G_1), \id_G)}(\lambda)$
$(\lambda\in L'(=G))$
\eqref{eqn:correspondence:dependdYB}.
\end{rem}

\end{document}